\documentclass[amssymb,11pt]{amsart}
\usepackage{amscd,amssymb,euscript,amsthm}
\theoremstyle{plain}
\newtheorem{thm}{Theorem}[section] 

\newtheorem{prop}[thm]{Proposition}
\newtheorem{lem}[thm]{Lemma}
\theoremstyle{definition}
\newtheorem{defn}[thm]{Definition}
\theoremstyle{remark}
\newtheorem{rem}[thm]{Remark}

\numberwithin{equation}{section}
\newcommand{\ran}{\operatorname{ran}}
\newcommand{\rank}{\operatorname{rank}}

\newcommand{\tr}{\operatorname{trace}}
\newcommand{\ind}{\operatorname{index}}
\newcommand{\co}{\operatorname{conv}}
\newcommand{\lm}{\operatorname{Lim{^1(X)}}}

\newcommand{\lum}{\operatorname{Lim{^1(\{0,1\})}}}
\def\<{\left<}
\def\>{\right>}
\def\cstar{$C^*$-algebra}
\begin{document}
\title[Diagonals of Normal Operators]{Diagonals of Normal Operators \\ with finite spectrum}
\author{William Arveson}
%
%
\address{Department of Mathematics,
University of California, Berkeley, CA 94720}
\email{arveson@math.berkeley.edu}

\begin{abstract}
Let $X=\{\lambda_1,\dots,\lambda_N\}$ be 
a finite set of complex numbers and 
let $A$ be a normal operator with spectrum $X$ that acts on 
a separable Hilbert space $H$.  
Relative to a fixed orthonormal basis $e_1,e_2,\dots$ for $H$, 
$A$ gives rise to a matrix 
whose diagonal is a sequence $d=(d_1,d_2,\dots)$ 
with the property that each of its terms $d_n$ belongs 
to the convex hull of $X$.  Not all 
sequences with that property can arise 
as the diagonal of a normal operator with spectrum $X$.  

The case where $X$ is a set of real numbers has 
received a great deal of attention 
over the years, and is reasonably well
(though incompletely) understood.    
In this paper we take up the case in which $X$ is the set of vertices of 
a convex polygon in $\mathbb C$.  The critical sequences $d$ 
turn out to be those  
that accumulate rapidly in $X$ in the sense that 
$$
\sum_{n=1}^\infty {\rm{dist}}\,(d_n,X)<\infty.   
$$ 
We show that there is an abelian group $\Gamma_X$ -- 
a quotient of $\mathbb R^2$ by a countable subgroup  
with concrete arithmetic properties -- 
and a surjective mapping of such sequences 
$d\mapsto s(d)\in\Gamma_X$ with the following 
property: If $s(d)\neq 0$, then $d$ is not the diagonal 
of any such operator $A$.  

We also show that while this is the only obstruction 
when $N=2$, there are other (as yet unknown) obstructions 
when $N=3$.  
\end{abstract}
\maketitle

\section{Introduction.}\label{S:pr}

Given a self-adjoint $n\times n$ matrix $A$, the diagonal of $A$ 
and the eigenvalue list of $A$ are two points of $\mathbb R^n$ that 
bear some relation to each other.  
The Schur-Horn theorem characterizes that relation in terms of 
a system of linear inequalities  
\cite{schurDiag},  \cite{hornDiag}.  That characterization has attracted 
a great deal of interest over the years, and has been 
generalized in remarkable ways.  For example, 
\cite{KostantSHthm}, \cite{atiyHorn},  \cite{GSI}, 
\cite{GSII} represent some of the milestones.  
More recently, a characterization of the diagonals 
of projections acting on infinite dimensional Hilbert spaces has been 
discovered \cite{kadPnasI},\cite{kadPnasII}, and a version of the 
Schur-Horn theorem for positive trace-class operators 
was given in \cite{arvKad1}.  
The latter reference contains a somewhat more complete 
historical discussion.  

Let $X$ be a finite subset of the complex plane $\mathbb C$, 
and consider the set $\mathcal N(X)$ of all normal operators 
acting on a separable Hilbert space $H$ that have spectrum 
$X$ with uniformly infinite multiplicity, 
$$
\mathcal N(X)=\{A\in\mathcal B(H): A^*A=AA^*,\quad \sigma(A)=\sigma_e(A)=X\}.  
$$
The set $\mathcal N(X)$ is invariant under the action of 
the group of $*$-automorphisms 
of $\mathcal B(H)$, and it is 
closed in the operator norm.  
Fixing an orthonormal basis $e_1,e_2,\dots$ for $H$, one may 
consider the (non-closed) set $\mathcal D(X)$ of all diagonals of 
operators in $\mathcal N(X)$ 
$$
\mathcal D(X)=\{(\langle Ae_1,e_1\rangle,\langle Ae_2,e_2\rangle,\dots))\in\ell^\infty: 
 A\in\mathcal N(X)\}.  
$$
In this paper we address the problem of determining the elements of $\mathcal D(X)$.  

Notice that for every sequence $d=(d_1,d_2,\dots)$ in $\mathcal D(X)$,  
{\em each term $d_n$ must belong to the convex hull of $X$.}  Indeed, since 
there is a normal operator $A$ with spectrum $X$ such that 
$d_n=\langle Ae_n,e_n\rangle$, $n\geq 1$,  
each $d_n$ must belong to the numerical range of $A$, and 
the closure of the numerical range of a normal operator is the convex hull 
of its spectrum.

This necessary condition $d_n\in\co X$, $n\geq 1$, is  
not sufficient.  Indeed, a 
characterization of $\mathcal D(\{0,1\})$ (the set of diagonals of projections) 
was given in \cite{kadPnasII}, the main 
assertion of which can be paraphrased 
as follows:

\begin{thm}[Theorem 15 of \cite{kadPnasII}]\label{prThm1}
Let $d=(d_1,d_2,\dots)\in\ell^\infty$ be a sequence 
satisfying $0\leq d_n\leq 1$ for every $n$ and 
$$
\sum_{n=1}^\infty d_n=\sum_{n=1}^\infty 1-d_n=\infty.
$$
Let $a,b\in[0,\infty]$ be the numbers 
$$
a=\sum_{d_n\leq 1/2} d_n,\qquad b=\sum_{d_n>1/2}1-d_n.
$$  
Then one has the following dichotomy: 
\begin{enumerate}
\item[(i)]If $a+b=\infty$, then $d\in\mathcal D(\{0,1\})$.  
\item[(ii)]If $a+b<\infty$, 
then $d\in\mathcal D(\{0,1\}) \iff a-b\in\mathbb Z$.  
\end{enumerate}
\end{thm}

In a recent paper\cite{arvKad1}, a 
related spectral characterization was found for 
the possible diagonals of positive trace-class operators.  That 
paper did not address the case of more general self-adjoint 
operators, and in particular, the results of \cite{arvKad1} 
shed no light on the phenomenon (ii) of Theorem \ref{prThm1}.  
This paper grew out of an effort to understand 
that phenomenon as an index obstruction.  We 
achieve that for certain finite subsets $X\subseteq \mathbb C$ 
in place of $\{0,1\}$, namely those that are affinely independent 
in the sense that none of the points of $X$ can be written as 
a nontrivial convex combination of the others - i.e., 
when $X$ is the set of vertices of a convex 
polygon.    

The basic issues taken up here bear some 
relation to A. Neumann's work on the infinite-dimensional 
Schur-Horn theorem for self-adjoint operators.  But there is a 
fundamental difference in the 
nature of the characterizations of \cite{aNeumann1} and the results 
below that goes beyond the fact that Neumann confines attention to 
self-adjoint operators.  The comparison is clearly seen 
for the two-point set $X=\{0,1\}$.  
In that case, the results of \cite{aNeumann1} provide the following description of 
the {\em closure} of $\mathcal D(X)$ in the $\ell^\infty$-norm:
$$
\overline{\mathcal D(X)}=\{d=(d_n)\in\ell^\infty: 0\leq d_n\leq 1,\quad 
n=1,2,\dots\},
$$ 
see  Lemma 2.13 and Proposition 3.12 of \cite{aNeumann1}.  
Thus, 
the exceptional cases described in part 
(ii) of Theorem \ref{prThm1} disappear 
when one passes from $\mathcal D(X)$ to its closure in the 
$\ell^\infty$-norm.  
In more explicit terms,  
while sequences $d=(d_n)$ satisfying $0\leq d_n\leq 1$, $n\geq 1$, 
$a+b<\infty$, $a-b\notin \mathbb Z$,  
fail to belong to $\mathcal D(X)$, they are all absorbed 
into its norm-closure.  

It is these ``exceptional" cases that 
we seek to understand here, for more 
general finite sets $X\subseteq \mathbb C$.  
Our main result (Theorem \ref{dgThm1} below) 
identifies an index obstruction corresponding to 
(ii) above when $X$ is the set of vertices of a convex 
polygon $P$.  Specifically, consider the set of all sequences 
$d=(d_n)$ that satisfy $d_n\in P$, $n=1,2,\dots$, and which 
accumulate rapidly in $X$ in the precise sense that 
$$
\sum_{n=1}^\infty {\rm{dist\,}}(d_n, X)<\infty.   
$$
We show that there is 
a discrete abelian group $\Gamma_X$, depending only on 
the arithmetic properties of $X$, and a surjective mapping 
$d\mapsto s(d)\in \Gamma_X$ of the set of all such 
sequences $d$, with the following property: If $s(d)\neq 0$, 
then $d$ is not the 
diagonal of any operator in $\mathcal N(X)$.  
We use Theorem \ref{prThm1} to show that 
this is the only obstruction in the case of two point sets;  
but we also show by example that there are other (as yet unkown) obstructions 
in the case of three-point sets.

Finally, I want to thank Richard Kadison, whose work 
\cite{kadPnasI}, \cite{kadPnasII} initially inspired 
this effort, and with whom I have had the pleasure of many 
helpful conversations.

\section{Sequences in $\lm$ and the group $\Gamma_X$}\label{S:ob}

Let $X=\{\lambda_1,\dots,\lambda_N\}$ be a finite set of complex numbers.  
For every complex number $z$ we write 
$$
d(z,X)=\min_{\lambda\in X}|z-\lambda|
$$
for the distance from $z$ to $X$.  
We consider the space $\lm$ of all sequences 
$a=(a_1, a_2, \dots)\in\ell^\infty$ with the property 
\begin{equation}\label{obEq1}
\sum_{n=1}^\infty d(a_n,X)<\infty.    
\end{equation}
Thus, a sequence $a=(a_n)$ belongs to $\lm$ iff all of its 
limit points belong to $X$ and it converges 
rapidly to its limit points in the following sense: 
there is a sequence $x=(x_n)$ satisfying 
$x_n\in X$ for every $n=1,2,\dots$, and 
\begin{equation}\label{obEq2}
\sum_{n=1}^\infty |a_n-x_n|<\infty.   
\end{equation}
In the context of Theorem \ref{prThm1}, (\ref{obEq2}) 
reduces to the hypothesis of (ii) when $X=\{0,1\}$, see 
Section \ref{S:cv}.  
In this section we show that every element $a\in\lm$ has 
a ``renormalized" sum that takes values in an 
abelian group $\Gamma_X$ naturally associated with $X$.

For fixed $a\in\lm$ there are many $X$-valued sequences 
$x=(x_n)$ that satisfy (\ref{obEq2}).   
Nevertheless, one can attempt to define a 
``renormalized" sum of an element $a\in\lm$ 
by choosing a sequence $x_n\in X$ 
that satisfies (\ref{obEq2}) and forming the complex 
number 
$$
s=\sum_{n=1}^\infty a_n-x_n.  
$$
While the value of $s$ depends on the choice 
of $x\in X$, the following observation shows that 
the ambiguity is associated 
with a countable subgroup of the additive group 
of $\mathbb C$.  

\begin{prop}\label{obProp1}Let $X=\{\lambda_1,\dots,\lambda_N\}$ 
be a finite subset of $\mathbb C$ and 
fix $a=(a_n)\in\lm$.  For any two sequences $x=(x_n),y=(y_n)$ of points 
in $X$ that satisfy 
$$
\sum_{n=1}^\infty |a_n-x_n|<\infty,\qquad \sum_{n=1}^\infty |a_n -y_n|<\infty,
$$ 
the sequence of differences $x-y=(x_n-y_n)$ is finitely nonzero,  and 
there are integers $\nu_1,\dots,\nu_N\in\mathbb Z$ 
such that  $\nu_1+\nu_2+\cdots+\nu_N=0$ and 
\begin{equation}\label{obEq3}
\sum_{n=1}^\infty x_n-y_n=\nu_1\lambda_1+\nu_2\lambda_2+\cdots+\nu_N\lambda_N.  
\end{equation}
\end{prop}

\begin{proof}Since $a-y$ and $a-x$ both belong to $\ell^1$, their 
difference $x-y$ must also belong to $\ell^1$.  
Since 
$x_n-y_n$ takes 
values in the finite set of differences $X-X$ and belongs to $\ell^1$, 
it must vanish for all but finitely many $n$, and for each of the remaining 
$n$, $x_n-y_n$ is of the form $\lambda_{i_n}-\lambda_{j_n}$ where 
$i_n,j_n\in\{1,2,\dots,N\}$.  It follows that  
$$
\sum_{n=1}^\infty x_n-y_n
$$
is a finite sum of terms of the form $\lambda_i-\lambda_j$, 
$1\leq i,j\leq N$, and such a number has the 
form (\ref{obEq3}) with integer coefficients $\nu_k\in\mathbb Z$ 
having sum $0$.  
\end{proof}

\begin{defn}[The obstruction group $\Gamma_X$]\label{obDef1}For every finite set 
of $N\geq 2$ complex numbers $X=\{\lambda_1,\dots,\lambda_N\}$, 
let $K_X$ be the additive subgroup of $\mathbb C$ consisting 
of all $z$ of the form 
$$
z=\nu_1\lambda_1+\cdots+\nu_N\lambda_N
$$
where $\nu_1,\dots,\nu_N\in\mathbb Z$ satisfy 
$\nu_1+\cdots+\nu_N=0$.   $\Gamma_X$ will denote the 
quotient of abelian groups 
$$
\Gamma_X=\mathbb C/K_X.
$$
\end{defn}

$K_X$ is the subgroup of $\mathbb C$ generated 
by the set of differences $\lambda_i-\lambda_j$, for $i, j=1,\dots,N$,
or equivalently by $\{\lambda_2-\lambda_1,\lambda_3-\lambda_1,\dots,\lambda_N-\lambda_1\}$.  
Hence the rank of $K_X$ is at most $N-1$.  The rank is $N-1$ iff 
when one views $\mathbb C$ as a vector space over the field $\mathbb Q$ of 
rational real numbers, the set of differences 
$\{\lambda_2-\lambda_1,\lambda_3-\lambda_1,\dots,\lambda_N-\lambda_1\}$
becomes a linearly independent set.  

By Proposition \ref{obProp1}, we can define a 
map $s: \lm\to \Gamma_X$ as follows:  For every 
$a=(a_n)\in\lm$, choose a sequence $x=(x_n)$ that takes values in 
$X$ and satisfies $\sum_n|a_n-x_n|<\infty$, and let   
$s(a)$ be the coset  
$$
s(a)=\sum_{n=1}^\infty a_n-x_n + K_X \in\Gamma_X.  
$$
\begin{defn}[Renormalized Sum]\label{obDef2}
For every sequence $a\in\lm$, the element 
$s(a)\in\Gamma_X$ is called the {\em renormalized sum}
of $a$.  
\end{defn}
\noindent
When it is necessary to call attention to the 
set  $X$ of vertices, we will write $s_X(a)$ rather than 
$s(a)$.

\begin{rem}[Surjectivity of the map $s: \lm\to \Gamma_X$]  
One thinks of $\Gamma_X$ as an uncountable discrete abelian 
group.  
It is easy to see that the map $s$ is surjective.  Indeed, for any 
$z\in\mathbb C$ the coset $z+K_X\in\Gamma_X$ is realized as 
the value $s(a)$ of a renormalized sum  
as follows.  Choose any sequence $u=(u_n)$ in $\ell^1$ such that 
$$
\sum_{n=1}^\infty u_n=z,    
$$
and let $x=(x_n)$ be an arbitrary sequence satisfying $x_n\in X$ for 
every $n\geq 1$.  Then the sum $a=u+x$  belongs to $\lm$ and 
satisfies $s(a)=z+K_X$.   
\end{rem}

We require the following elementary description of 
sequences in $\lm$.  

\begin{prop}\label{obProp2}
Let $X=\{\lambda_1,\dots,\lambda_N\}$ be a finite set of 
complex numbers and let $a=(a_n)\in\ell^\infty$.  Then the following 
are equivalent:
\begin{enumerate}
\item[(i)]$a\in\lm$.  
\item[(ii)]One has the summability condition
$$
\sum_{n=1}^\infty |f(a_n)|<\infty, 
$$
where $f$ is the polynomial 
$f(z)=(z-\lambda_1)(z-\lambda_2)\cdots(z-\lambda_N)$.  
\end{enumerate}
\end{prop}

\begin{proof}
Let 
$$
\delta=\min_{1\leq i< j\leq N}|\lambda_i-\lambda_j|  
$$
be the minimum distance between distinct points of $X$.  
Note first that whenever $z\in\mathbb C$ satisfies 
$d(z,X)\leq \delta/2$, one has 
\begin{equation}\label{obEq4}
|f(z)|\geq d(z,X)(\delta/2)^{N-1}.  
\end{equation}
Indeed, if we choose $k$ such that $d(z,X)=|z-\lambda_k|$, then
for $j\neq k$ we have 
$$
|z-\lambda_j|\geq |\lambda_k-\lambda_j|-|z-\lambda_k|\geq 
\delta-\delta/2=\delta/2, 
$$
hence 
$$
|(z-\lambda_1)\cdots(z-\lambda_N)|\geq |z-\lambda_k|(\delta/2)^{N-1},
$$
and (\ref{obEq4}) follows.  

If $a=(a_1,a_2,\dots)$ is a sequence such that $\sum_n|f(a_n)|$ converges, 
then since $|f(z)|\geq d(z,X)^N$ for all $z\in\mathbb C$, it follows that 
$d(a_n,X)^N\to 0$ as $n\to\infty$, hence there is an $n_0$ such that 
$d(a_n,X)\leq \delta/2$ for $n\geq n_0$.  (\ref{obEq4}) implies 
$$
\sum_{n=n_0}^\infty |f(a_n)|\geq (\delta/2)^{N-1}\sum_{n=n_0}^\infty d(a_n,X), 
$$
so that $\sum_n d(a_n,X)$ converges, hence $a\in\lm$.  

Conversely, assuming that $\sum_nd(a_n,X)$ converges, let 
$$
R=\max(|\lambda_1|,\dots,|\lambda_N|)>0.
$$  
Since $d(a_n,X)\to 0$, we can find $n_1$ such that 
$|a_n|\leq 2R$ for $n\geq n_1$.  Choosing $k_1, k_2,\dots$ so 
that $d(a_n,X)=|a_n-\lambda_{k_n}|$ we have 
$$
|f(a_n)|=|a_n-\lambda_1|\cdots|a_n-\lambda_N|\leq 
|a_n-\lambda_{k_n}|(3R)^{N-1}=d(a_n,X)(3R)^{N-1}
$$
for all $n\geq n_1$, hence $\sum_n|f(a_n)|$ converges.  
\end{proof}

\section{$X$-decompositions}\label{S:xd}

We view  
$\ell^\infty$ as a commutative \cstar\ with unit $\mathbf 1$, 
and the elements 
of $\ell^\infty$ as bounded functions $a: \mathbb N\to \mathbb C$, 
with norm 
$$
\|a\|=\sup_{n\geq 1}|a(n)|.  
$$
Let $P$ be a convex polygon in the complex plane 
and let 
$X=\{\lambda_1,\dots,\lambda_N\}$ be 
its set of vertices, with $N\geq 2$.  
In this section we describe an elementary (nonunique) decomposition for 
sequences that take values in $P$, and we show that under 
certain circumstances, all such decompositions of a 
sequence must share a key property.

\begin{prop}\label{xdProp1}Let $P\subseteq \mathbb C$ be a convex polygon with 
vertices $\{\lambda_1,\dots,\lambda_N\}$.  
Every sequence $a\in\ell^\infty$ satisfying 
$a(n)\in P$,  $n\geq 1$, can be decomposed into a sum of the form
\begin{equation}\label{xdEq2}
a=\lambda_1e_1+\lambda_1e_2+\cdots+\lambda_Ne_N
\end{equation}
where $e_1,\dots,e_N$ are positive elements of $\ell^\infty$ 
satisfying $e_1+\cdots+e_N=\mathbf 1$.  

Conversely, any sequence $a$ of the form 
(\ref{xdEq2}), with positive elements $e_k$ summing to $\mathbf 1$,  
must satisfy $a(n)\in P$ for every $n\geq 1$.  
\end{prop}

\begin{proof}
Fix $a$ and choose $n\geq 1$.  Since $a(n)\in P$ and 
$P$ is the convex hull of $\{\lambda_1,\dots,\lambda_N\}$, we can find 
a point $(e_1(n),\dots,e_N(n))\in\mathbb R^N$ such 
that $e_1(n)\geq 0,\dots,e_N(n)\geq 0$, $e_1(n)+\cdots+e_N(n)=1$, and 
$$
a(n) = \sum_{k=1}^N e_k(n)\lambda_k.
$$
The sequences $e_k=(e_k(1),e_k(2),\dots)\in\ell^\infty$
satisfy $e_k\geq 0$, $e_1+\cdots+e_N=\mathbf 1$, 
and the asserted representation (\ref{xdEq2}) follows.  

The converse assertion is obvious.  
\end{proof}

\begin{defn}\label{xdDef1}Let $P\subseteq \mathbb C$ be a convex 
polygon whose set of vertices is $X=\{\lambda_1,\dots,\lambda_N\}$, 
and let $a\in\ell^\infty$ satisfy $a(n)\in P$, $n\geq 1$.  
A representation of the form (\ref{xdEq2}) is called 
an {\em $X$-decomposition} of $a$.  
\end{defn}

Despite the fact that $X$-decompositions are not unique 
except in very special circumstances, 
there is a property common to 
all $X$-decompositions of $a$ in cases where 
$$
\sum_{n=1}^\infty d(a(n),X)<\infty.  
$$
That result (Theorem \ref{xdThm1} below) requires the case $n=2$ of the following: 

\begin{lem}\label{xdLem1}
Let $P$ be a convex polyhedron in $\mathbb R^n$ with 
extreme points $x_1,\dots,x_r$.  Consider the simplex 
$$
\Delta=\{(t_1,\dots,t_r)\in\mathbb R^r: t_i\geq 0,\ t_1+\cdots+t_r=1\}
$$
and the affine map of $\Delta$ onto $P$ defined by 
$$
t\in\Delta\mapsto x(t)=t_1x_1+\cdots+t_rx_r.  
$$
For any choice of norms on $\mathbb R^r$ and $\mathbb R^n$, there 
is a constant $C>0$ such that 
\begin{equation}\label{xdEq2.1}
d(t,\{\delta_1,\dots,\delta_r\})\leq C\cdot d(x(t),\{x_1,\dots,x_r\}), 
\qquad t\in\Delta, 
\end{equation}
where $d(v,S)=\inf \{\|v-s\|: s\in S\}$ denotes the distance from 
a vector $v$ to a set $S$, and where $\delta_1,\dots, \delta_r$ are 
the extreme points of $\Delta$, $(\delta_k)_j=\delta_{kj}$.  
\end{lem}

\begin{proof} 
It suffices to show that for each $k=1,\dots, r$, 
there is a constant $C_k$ such that 
\begin{equation}\label{xdEq3}
\|t-\delta_k\|\leq C_k\cdot\|x(t)-x_k\|,\qquad t\in\Delta.  
\end{equation}
Indeed, (\ref{xdEq3}) implies 
$$
\min_{1\leq k\leq r}\|t-\delta_k\|\leq C\min_{1\leq k\leq r}\|x(t)-x_k\|,
\qquad t\in\Delta,
$$
where $C=\max(C_1,\dots,C_r)$, and (\ref{xdEq2.1}) follows.  

By symmetry, it suffices to prove 
(\ref{xdEq3}) for $k=1$; moreover, after performing an affine 
translation if necessary, there is 
no loss of generality if we assume that $x_1=0$ is one of 
the extreme points of $P$.  
For every extreme point $e$ of a convex polyhedron $P$, there is a supporting 
hyperplane that meets $P$ only at $\{e\}$.  
Thus there 
is a linear functional $f$ on $\mathbb R^n$ such that 
$f(x)>0$ for all nonzero $x\in P$.    For each $t\in\Delta$ we 
have 
$$
f(x(t))=\sum_{k=2}^rt_kf(x_k)\geq \min_{2\leq k\leq r} f(x_k)\cdot\sum_{k=2}^r t_k.  
$$
After noting that $f(x(t))\leq \|f\|\cdot \|x(t)\|$,  we obtain 
$$
\|t-\delta_k\|=\sum_{k=1}^r|t_k-\delta_1(k)|=2\sum_{k=2}^r t_k\leq 
\frac{2\cdot\|f\|}{\min(f(x_2),\dots,f(x_r))}\cdot \|x(t)\|, 
$$
and (\ref{xdEq3}) follows after noting that the left side 
of the preceding inequality 
dominates $\epsilon\cdot\|t-\delta_1\|$ for an appropriately small  
positive constant $\epsilon$.  
\end{proof}

\begin{thm}\label{xdThm1}Let $X=\{\lambda_1,\dots,\lambda_N\}$ be 
the set of vertices of a convex polygon $P\subseteq \mathbb C$ and 
let $a\in\ell^\infty$ satisfy $a(n)\in P$ for $n\geq 1$.  
If $a\in\lm$, then for every $X$-decomposition of the form (\ref{xdEq2})
$$
a=\lambda_1e_1+\cdots+\lambda_Ne_N, 
$$
each of the sequences $e_1,\dots,e_N$ 
belongs to ${\rm{Lim}}{^1\{0,1\}}$.  
\end{thm}

\begin{proof} Consider the Euclidean norm on $\mathbb C$, 
the norm 
$$
\|(x_1,\dots,x_N)\|=|x_1|+\cdots+|x_N|
$$ 
on $\mathbb R^N$, and fix $n=1,2,\dots$.   
Since the point $(e_1(n),\dots,e_N(n))\in\mathbb R^N$ belongs to the simplex 
$\Delta$ of Lemma \ref{xdLem1}, there is a constant 
$C>0$ such that 
\begin{align*}
d((e_1(n), \dots,e_N(n)),\{\delta_1,\dots,\delta_N\})
&\leq C \cdot d(a(n),\{\lambda_1,\dots,\lambda_N\})
\\&
=C\cdot d(a(n),X),
\qquad n\geq 1.  
\end{align*}
Since for every point $t=(t_1,\dots,t_N)$ in the simplex $\Delta$ and 
for every fixed $k=1,\dots, N$ we 
have 
\begin{align*}
d(t_k,\{0,1\})&=\min(t_k,1-t_k)
\\&
\leq 
\min_{1\leq k\leq N}(\sum_{j\neq k}t_j+(1-t_k))=d(t,\{\delta_1,\dots,\delta_N\}), 
\end{align*}
it follows that 
\begin{align*}
d(e_k(n),\{0,1\})&\leq d((e_1(n), \dots,e_N(n)), \{\delta_1,\dots,\delta_N\})
\\&
\leq 
C d(a(n),\{\lambda_1,\dots,\lambda_N\})
=Cd(a(n),X).
\end{align*}
Using $a\in\lm$, we can sum the preceding inequality on 
$n$ to obtain 
$$
\sum_{n=1}^\infty d(e_k(n),\{0,1\})\leq C\cdot \sum_{n=1}^\infty d(a(n),X)<\infty.  
$$
Hence $e_k\in{\rm{Lim}}{^1\{0,1\}}$.   
\end{proof}

\section{Two Projections.}\label{S:tp}

It is known that for any pair of projections $P,Q\in\mathcal B(H)$ 
for which $P-Q$ is a trace-class operator, 
$\tr(P-Q)$ must be an integer.  For example, 
the result can be found in Effros' article 
(\cite{EdIntell}, see Lemma 4.1).  
We require an appropriate extension of that 
result to the case where $P-Q$ is 
merely a Hilbert-Schmidt operator, Theorem \ref{tpThm1} below.  
Throughout the remainder 
of this paper, we write $\mathcal L^1$ (resp. $\mathcal L^2$) 
for the Banach space of trace-class operators (resp. Hilbert-Schmidt 
operators) acting on a given Hilbert space, and 
we write $P^\perp$ for $\mathbf 1-P$ when $P\in\mathcal B(H)$ is a projection.

\begin{thm}\label{tpThm1}
Let $M$, $N$ be subspaces of a Hilbert space $H$ with 
respective projections $P,Q$, and assume that $P-Q\in\mathcal L^2$.  

Then both $Q(P-Q)Q$ and 
$Q^\perp(P-Q)Q^\perp$ belong to $\mathcal L^1$, 
both subspaces $M\cap N^\perp$ and 
$N\cap M^\perp$ are finite-dimensional,  and  
\begin{equation}\label{tpEq1}
\tr (Q(P-Q)Q+Q^\perp(P-Q)Q^\perp)=\dim(M\cap N^\perp)-\dim(N\cap M^\perp).  
\end{equation}
In particular, $QPQ+Q^\perp PQ^\perp -Q$ is a trace-class operator 
such that 
\begin{equation}\label{tpEq1.1}
\tr(QPQ+Q^\perp PQ^\perp -Q)\in\mathbb Z.
\end{equation}
\end{thm}

In the proof we will show that the left side of (\ref{tpEq1}) is the 
index of a Fredholm operator, and for that we require the following 
elementary result for which we lack a convenient reference: 

\begin{lem}\label{tpLem1}
Let $H$, $K$ be Hilbert spaces and let $A:H\to K$  
be an operator such that both  $\mathbf 1_H-A^*A$ and 
$\mathbf 1_K-AA^*$ are trace-class.  Then $A$ is a Fredholm 
operator in $\mathcal B(H,K)$ whose index is given by the formula 
\begin{equation}\label{tpEq2}
\ind A=\tr(\mathbf 1_H-A^*A)-\tr(\mathbf 1_K-AA^*)
\end{equation}
\end{lem}

\begin{proof}
Consider the polar 
decomposition $A=UB$, where $B$ is a positive operator  and 
$U$ is a partial isometry with initial space 
$\ker A^\perp$ and range  $\overline{\ran A}=\ker A^{*\perp}$.  
Let $C$ be the restriction of $B^2$ to $\ker A^\perp$.   
Then $A^*A=B^2=C\oplus 0_{\ker A}$ and $AA^*=C^\prime\oplus 0_{\ker A^*}$, 
where $C^\prime$ is unitarily equivalent to $C$; indeed, the 
restriction of $U$ to $\ker A^\perp$ implements a unitary equivalence 
of $C$ and $C^\prime$.  

It follows that 
$
\mathbf 1-A^*A=(\mathbf 1-C)\oplus \mathbf 1_{\ker A}, 
$
and 
$\mathbf 1-AA^*$ is unitarily equivalent to $(\mathbf 1-C)\oplus \mathbf 1_{\ker A^*}$.  
Since $\mathbf 1-C$ is a 
trace-class operator, we have $\tr(\mathbf 1-A^*A)=\tr(\mathbf 1-C)+\dim\ker A$, 
and similarly $\tr(\mathbf 1-AA^*)=\tr(\mathbf 1-C)+\dim\ker A^*$.  
The terms involving $\tr(\mathbf 1-C)$ cancel, and 
\begin{equation*}
\tr(\mathbf 1-A^*A)-\tr(\mathbf 1-AA^*)=\dim\ker A-\dim \ker A^*
=\ind A, 
\end{equation*}
as asserted.  
\end{proof}

\begin{proof}[Proof of Theorem \ref{tpThm1}]  We claim first 
that $Q(P-Q)Q\in\mathcal L^1$.  Indeed, we have 
$Q(P-Q)Q=-(Q-QPQ)$, and $Q-QPQ$ is a positive operator satisfying 
$$
\tr (Q-QPQ)=\tr QP^\perp Q=\tr |P^\perp Q|^2=\tr |(P-Q)Q|^2<\infty.  
$$
Similarly, $Q^\perp(P-Q)Q^\perp =Q^\perp PQ^\perp$, and 
$$
\tr Q^\perp PQ^\perp=\tr|PQ^\perp|^2=\tr|(P-Q)Q^\perp|^2<\infty, 
$$
so that $Q^\perp(P-Q)Q^\perp\in\mathcal L^1$.  

Let $H_0$ be 
the subspace of $H$ spanned by the mutually orthogonal 
subspaces $M\cap N^\perp$ and $N\cap M^\perp$.  The restriction 
of $P-Q$ to $H_0$ is unitary with eigenvalues $\pm 1$, 
hence $\dim H_0=\tr|(P-Q)P_{H_0}|^2\leq \tr |P-Q|^2$ must be finite.  

Consider the operator $A: N\to M$ defined 
by restricting $P$ to $N=QH$.  
Obviously, $\ker A=N\cap M^\perp$ 
and $\ker A^*=M\cap N^\perp$, and we claim that $A$ satisfies 
the hypotheses of Lemma \ref{tpLem1}.  Indeed, 
\begin{align*}
(\mathbf 1_{N}-A^*A)Q&=Q-QPQ=QP^\perp Q=|P^\perp Q|^2\\
(\mathbf 1_{M}-AA^*) P&=P-PQP=PQ^\perp P =|Q^\perp P|^2.  
\end{align*}
Since $P-Q$ is Hilbert-Schmidt, 
$Q^\perp P=(P-Q)P$ and $P^\perp Q=(Q-P)Q$ are both 
Hilbert-Schmidt, hence $|Q^\perp P|^2$ and $|P^\perp Q|^2$ 
are both trace-class.  Thus we can apply (\ref{tpEq2}) 
to obtain 
\begin{align*}
\ind A&=\tr PQ^\perp P-\tr QP^\perp Q=
\tr Q^\perp P Q^\perp -\tr Q P^\perp Q
\\&
=\tr Q^\perp(P-Q)Q^\perp+\tr Q(P-Q)Q.  
\end{align*}
The left side is $\dim\ker A-\dim\ker A^*=\dim(N\cap M^\perp) -\dim(M\cap N^\perp)$, 
and (\ref{tpEq1}) follows from the preceding formula.  
\end{proof}

\section{Projections with diagonals in ${\rm{Lim}}{^1\{0,1\}}$.}\label{S:pl}

In this section we characterize the projections in $\mathcal B(H)$ 
whose diagonals relative to a given orthonormal basis belong to 
${\rm{Lim}}{^1\{0,1\}}$.

\begin{prop}\label{plProp1}
Let $e_1, e_2,\dots$ be an orthonormal basis for a Hilbert space 
$H$, let $\mathcal A$ be the maximal abelian von Neumann algebra 
of all operators that are diagonalized by  $(e_n)$, and 
let $E:\mathcal B(H)\to \mathcal A$ be the trace-preserving 
conditional expectation
\begin{equation}\label{plEq1}
E(A)=\sum_{n=1}^\infty \langle Ae_n,e_n\rangle e_n\otimes \bar{e_n}.  
\end{equation}

For every projection $P\in\mathcal B(H)$, the following are equivalent:
\begin{enumerate}
\item[(i)] The diagonal of $P$ relative to the basis $(e_n)$ 
belongs to ${\rm{Lim}}{^1\{0,1\}}$.  
\item[(ii)]$E(P)-E(P)^2\in\mathcal L^1.$
\item[(iii)]$P\in\mathcal A+\mathcal L^2$.  
\end{enumerate}
\end{prop}

The proof of Proposition \ref{plProp1} requires the following formula: 

\begin{lem}\label{plLem1}
Let $E:\mathcal B(H)\to \mathcal A$ be the map (\ref{plEq1}).  
Then for every projection $P\in\mathcal B(H)$ we have 
\begin{equation}\label{plEq2}
\tr\, (P-E(P))^2=\tr\,(E(P)-E(P)^2).  
\end{equation}
\end{lem}

\begin{proof}[Proof of Lemma \ref{plLem1}]
We have 
$E(P)e_n=d_ne_n$, where $d_n=\langle Pe_n,e_n\rangle$.  
Since $d_ne_n$ is the projection of $Pe_n$ onto 
the one-dimensional space $\mathbb C\cdot e_n$, we have 
$$
\|Pe_n-d_ne_n\|^2=\|Pe_n\|^2-\|d_ne_n\|^2=\langle Pe_n,e_n\rangle-d_n^2=d_n-d_n^2.  
$$
Hence 
$$
\tr\,[P-E(P)]^2=\sum_{n=1}^\infty \|Pe_n-d_ne_n\|^2
=\sum_{n=1}^\infty d_n-d_n^2, 
$$
and the right side is evidently the trace of  $E(P)-E(P)^2$.  
\end{proof}

\begin{proof}[Proof of Proposition \ref{plProp1}]
Let $d=(d_1,d_2,\dots)$ be the diagonal of $P$ relative to $(e_n)$, 
$d_n=\langle Pe_n,e_n\rangle$, $n\geq 1$.  Then 
$$
\tr(E(P)-E(P)^2)=\sum_{n=1}^\infty d_n-d_n^2, 
$$
hence the equivalence of (i) and (ii) follows from 
Proposition \ref{obProp2}.  

(iii)$\implies$(ii): Assume first that the projection $P$ can be decomposed into a sum 
$P=A+T$ where $A\in\mathcal A$ and $T$ is Hilbert-Schmidt.  Then 
$P-E(P)=T-E(T)$, and $T-E(T)$ is a Hilbert-Schmidt operator.  By 
(\ref{plEq2}), we obtain 
$$
\tr\,(E(P)-E(P)^2)=\tr\,(P-E(P))^2=\tr\,(T-E(T))^2<\infty.  
$$

(ii)$\implies$(iii): 
Assume  
that $\tr(E(P)-E(P)^2))<\infty$, and consider the operator 
$T=P-E(P)$.  By (\ref{plEq2}), 
we have 
$$
\tr\,T^2=\tr\,(E(P)-E(P)^2)<\infty, 
$$
so that $T$ is Hilbert-Schmidt.  Thus $P=E(P)+T\in\mathcal A+\mathcal L^2$.  
\end{proof}

\section{Diagonals of operators in $\mathcal N(X)$.}\label{S:dg}
We are now in position to prove our main result.  Let 
$X=\{\lambda_1,\dots,\lambda_N\}$ be the set of vertices 
of a convex polygon $P\subseteq\mathbb C$ and let $\mathcal N(X)$ 
be the set of all normal operators $A$ acting on a separable Hilbert space 
$H$ that have spectrum $X$ with infinite multiplicity
$$
\sigma(A)=\sigma_e(A)= X.
$$
Fix an orthonormal basis $e_1,e_2,\dots$ for $H$.  
There are two necessary conditions 
that a sequence $d=(d_1,d_2,\dots)\in\ell^\infty$ must satisfy 
for it to be the diagonal of an operator $A\in\mathcal N(X)$, 
$d_n=\langle Ae_n,e_n\rangle$, $n\geq 1$, namely: 
\begin{enumerate}
\item[(i)]$d_n\in P$ for every $n\geq 1$, 
\item[(ii)] $d$ has an $X$-decomposition 
$$
d=\lambda_1 E_1+\cdots+\lambda_N E_N
$$
in which $\sum_{n=1}^\infty E_k(n)=\infty$ for every $k=1,\dots,N$.  
\end{enumerate}
Indeed the projections $P_k$ arising from the spectral
representation of $A$ 
$$
A=\lambda_1P_1+\cdots+\lambda_N P_N
$$
have diagonals 
$E_k(n)=\langle P_ke_n,e_n\rangle$ 
that give rise to an $X$-decomposition 
with  the property (ii).  
The requirements (i), (ii) on a sequence do not 
guarantee that it is the diagonal of an operator in $\mathcal N(X)$.  
We now identify an  
obstruction that emerges when $d\in\lm$ 
and which involves the renormalized sum $s:\lm\to \Gamma_X$ of  
Definition \ref{obDef1}. .

\begin{thm}\label{dgThm1} Let $X=\{\lambda_1,\dots,\lambda_N\}$ 
be the set of vertices of a 
convex polygon $P\subseteq \mathbb C$ and let $d=(d_1,d_2,\dots)$ be 
a sequence of complex numbers satisfying $d_n\in P$, 
$n\geq 1$, together with the summability condition 
\begin{equation}\label{dgEq1}
\sum_{n=1}^\infty |f(d_n)|<\infty,  
\end{equation}
where 
$f(z)=(z-\lambda_1)(z-\lambda_2)\cdots(z-\lambda_N)$.  
Then $d\in\lm$; and if   
$d$ is the diagonal of an 
operator in $\mathcal N(X)$, then $s(d)=0$.  
\end{thm}

\begin{proof} By Proposition \ref{obProp2}, the summability 
condition (\ref{dgEq1}) characterizes sequences in $\lm$.  

Fix an orthonormal basis $e_1,e_2,\dots$ for a Hilbert 
space $H$, and 
assume that there is an operator $A\in\mathcal N(X)$ 
such that $d_n=\langle Ae_n,e_n\rangle$, $n=1,2,\dots$.  In order to show 
that $s(d)=0$, we must find a sequence $(b_n)$ that takes values 
in $X$, satisfies $\sum_n|d_n-b_n|<\infty$, 
and we  must exhibit integers $\nu_1,\dots,\nu_n$ satisfying $\nu_1+\cdots+\nu_N=0$, 
and 
$$
\sum_{n=1}^\infty d_n-b_n=\nu_1\lambda_1+\cdots+\nu_N\lambda_N.
$$
For that, consider the maximal abelian algebra 
$\mathcal A$ of all operators that are diagonalized by the  
basis $e_1,e_2,\dots$, let $E:\mathcal B(H)\to\mathcal A$ be the trace-preserving 
conditional expectation
$$
E(T)=\sum_{n=1}^\infty \langle Te_n,e_n\rangle e_n\otimes\overline{e_n}, 
$$
and let $D=E(A)\in\mathcal A$ be the operator 
$$
D=\sum_{n=1}^\infty d_n\, e_n\otimes \overline{e_n}=
\sum_{n=1}^\infty \langle Ae_n,e_n\rangle  e_n\otimes \overline{e_n}.  
$$
We must find an operator $B\in\mathcal A\cap \mathcal N(X)$ such 
that $D-B\in\mathcal L^1$, and integers $\nu_1,\dots,\nu_N$ summing 
to zero, such that 
\begin{equation}\label{dgEq2}
\tr (D-B)=\nu_1\lambda_1+\cdots+\nu_N\lambda_N.
\end{equation}

The latter are achieved as follows.  Let 
$$
A=\lambda_1 P_1+\cdots+\lambda_N P_N
$$
be the spectral representation of $A$, with $P_1,\dots,P_N$ a set 
of mutually orthogonal infinite rank projections with sum 
$\mathbf 1$.  Then we have 
$$
D=E(A)=\lambda_1E(P_1)+\cdots+\lambda_N E(P_N), 
$$
so that $E(P_1),\dots,E(P_N)$ define an $X$-decomposition 
of $D$. Theorem \ref{xdThm1} implies that when one views 
the operators $E(P_k)$ as sequences in $\ell^\infty$, one has 
$E(P_k)\in\lum$ for each $k=1,\dots,N$.  

We claim that there is a sequence $Q_1,\dots,Q_N$ of mutually 
orthogonal projections in $\mathcal A$ having sum $\mathbf 1$ 
which satisfy 
\begin{equation}\label{dgEq3}
E(P_k)-Q_k\in\mathcal L^1,\qquad k=1,\dots,N.  
\end{equation}
Indeed, since $E(P_k)\in\lum$ for each $k$, the definition 
of $\lum$ implies that we can 
find projections $Q_1^0,\dots,Q_N^0\in\mathcal A$ such 
that 
$$
E(P_k)-Q_k^0\in\mathcal L^1,\qquad k=1,\dots,N.  
$$  
Considering each $Q_k^0$ as a sequence in $\ell^\infty$ that 
takes values in $\{0,1\}$, the sum $Q_1^0+\cdots+Q_N^0$ is 
a sequence taking values in $\{0,1,2,\dots,N\}$.  Consider 
the set  
$S=\{n\in\mathbb N: \sum_{k=1}^N Q_k^0(n)=1\}\subseteq \mathbb N$.  
Since $P_1+\cdots+P_N=\mathbf 1$ 
we have $E(P_1)+\cdots+E(P_N)=\mathbf 1$, and hence 
$$
\mathbf 1-\sum_{k=1}^N Q_k^0=\sum_{k=1}^N(E(P_k)-Q_k^0)\in\mathcal L^1.  
$$
It follows that 
$$
\sum_{n\notin S}|1-\sum_{k=1}^NQ_k^0(n)|<\infty.    
$$
The latter implies that $\mathbb N\setminus S$ is 
a finite set, and that $Q_1^0\cdot\chi_S,\dots,Q_N^0\cdot\chi_S$ 
are mutually orthogonal projections with sum $\chi_S$.  
Thus if we modify the sequence $Q_1^0,\dots,Q_N^0$ as follows, 
$$
Q_1=Q_1^0\cdot\chi_S+\chi_{\mathbb N\setminus S},
\ Q_2=Q_2^0\cdot\chi_S,\,\dots,\ Q_N=Q_N^0\cdot\chi_S,
$$ 
we obtain a new sequence of projections $Q_1,\dots,Q_N\in\mathcal A$ which 
are mutually orthogonal, have sum $\mathbf 1$, and satisfy 
(\ref{dgEq3}).  

Note too that since $\tr P_k=\rank P_k=\infty$ for every 
$k$, (\ref{dgEq3}) implies that $\tr Q_k=\rank Q_k=\infty$ as well.  
It follows that the operator 
$$
B=\lambda_1Q_1+\cdots+\lambda_NQ_N
$$
belongs to $\mathcal A$, satisfies $\sigma(B)=\sigma_e(B)=X$, and by  
construction, 
\begin{equation}\label{dgEq4}
D-B=E(A)-B\in\mathcal L^1.  
\end{equation}

It remains to show that $\tr(D-B)$ satisfies (\ref{dgEq2}) for 
integers $\nu_k$ as described there.  Indeed, since 
$D-B=\sum_k\lambda_k(E(P_k)-Q_k)$ and 
$E(P_k)-Q_k$ belongs to $\mathcal L^1$, we have  
\begin{equation*}
\tr(D-B)=
\sum_{k=1}^N\lambda_k\cdot \tr(E(P_k)-Q_k), 
\end{equation*}
so it suffices to show that 
\begin{equation}\label{dgEq5}
\tr(E(P_k)-Q_k)\in\mathbb Z,\qquad k=1,\dots,N, 
\end{equation}
and that the sum of  the $N$ integers of (\ref{dgEq5}) is $0$.  

In order to prove (\ref{dgEq5}) we appeal to Theorem \ref{tpThm1}.  
Note first that $P-E(P)$ belongs to $\mathcal L^2$.  Indeed, since 
$E(P_k)-Q_k$ is trace-class and $Q_k$ is a projection, we have 
$\tr (E(P)-E(P)^2)<\infty$, so by (\ref{plEq2}), 
$$
\tr[P-E(P)]^2=\tr(E(P)-E(P)^2)<\infty.  
$$
Since $E(P_k)-Q_k\in\mathcal L^1\subseteq \mathcal L^2$, we 
obtain 
$$
P_k-Q_k=(P_k-E(P_k))+E(P_k)-Q_k\in\mathcal L^2, \qquad k=1,\dots,N.
$$ 
From Theorem 
\ref{tpThm1} we conclude that 
$$
Q_kP_kQ_k+Q_k^\perp P_kQ_k^\perp -Q_k\in \mathcal L^1, 
$$
and moreover 
$$
\nu_k=\tr(Q_kP_kQ_k+Q_k^\perp P_kQ_k^\perp -Q_k)\in\mathbb Z.  
$$
Since $E(\mathcal L^1)\subseteq\mathcal L^1\cap\mathcal A$ and since 
\begin{align*}
E(Q_kP_kQ_k+Q_k^\perp P_kQ_k^\perp -Q_k)&=
Q_kE(P_k)Q_k+Q_k^\perp E(P_k)Q_k^\perp -Q_k
\\
&=E(P_k)-Q_k, 
\end{align*}
we find that $E(P_k)-Q_k\in\mathcal L^1$ and 
$$
\tr (E(P_k)-Q_k)=\tr(Q_kP_kQ_k+Q_k^\perp P_kQ_k^\perp -Q_k)=\nu_k\in\mathbb Z.  
$$
Since $\sum_{k=1}^N(E(P_k)-Q_k)=E(\mathbf 1)-\mathbf 1=0$, 
we have $\nu_1+\cdots+\nu_N=0$.  Finally, 
$$
\tr(D-B)=\sum_{k=1}^N\lambda_k\cdot\tr(E(P_k)-Q_k)=
\sum_{k=1}^N\lambda_k\nu_k, 
$$
and (\ref{dgEq2}) follows.  
\end{proof}

\section{Concluding Remarks, and an example}\label{S:cv}

We point out that Theorem \ref{dgThm1} specializes  
to the assertion (ii)  $\Rightarrow$ of 
Theorem \ref{prThm1} in the case $X=\{0,1\}$.  
Indeed, a straightforward calculation shows that 
for the two-point set $X=\{0,1\}$ one has 
$K_{\{0,1\}}=\mathbb Z$, so that 
$\Gamma_{\{0,1\}}=\mathbb C/\mathbb Z=\mathbb T\times \mathbb R$.  
Now the hypothesis of Theorem \ref{prThm1} (ii)
is that $a+b<\infty$, where 
$a$ and $b$ are defined by 
$$
a=\sum_{d_n\leq 1/2} d_n,\qquad b=\sum_{d_n>1/2}1-d_n.   
$$  
Let $(x_n)$ be the sequence 
$x_n=0$ when $d_n\leq 1/2$ and $x_n=1$ when 
$d_n>1/2$.  Then $a+b$ is finite $\iff$
$\sum_n|d_n-x_n|<\infty$ $\iff$ $d\in\lum$.  Moreover, 
$a-b=\sum_n d_n-x_n$, so that $a-b\in\mathbb Z$ $\iff$ $s_{\{0,1\}}(a)=0$.  

For more general sets $X$, 
the converse of Theorem \ref{dgThm1} would assert: 
\vskip0.05in

{\em Let $X$ be the set of vertices of a convex 
polygon and 
let $d$ be a sequence in $\lm$ such that 
$s_X(d)=0$.  Then there is an operator $N\in\mathcal N(X)$ and 
an orthonormal basis $e_1, e_2,\dots$ for $H$ such that 
$d_n=\langle Ne_n,e_n\rangle$, $n\geq 1$.  
}
\vskip0.05in

We first point out  that this converse is true 
when $X$ consists of just two points.  
To sketch the argument briefly, suppose 
$X=\{\lambda_1,\lambda_2\}$ with $\lambda_1\neq \lambda_2$.  
One can find an affine bijection $z\mapsto az+b$ of $\mathbb C$ 
that carries $\lambda_1$ to 
$0$ and $\lambda_2$ to $1$, and which therefore carries sequences 
in $\lm$ to sequences in $\lum$.  After noting that the 
operator mapping 
$T\mapsto aT+b\mathbf 1$ carries $\mathcal N(X)$ to 
$\mathcal N(\{0,1\})$, one can make use of 
the implication $\Leftarrow$ of Theorem \ref{prThm1} (ii) 
in a straightforward way to deduce the required result.  

On the other hand, the following example shows that this 
converse of Theorem \ref{dgThm1} fails for three-point 
sets.  

\begin{prop}\label{cvProp1}
Let $X=\{0,1,i\}$, $i$ denoting $\sqrt{-1}$,  and consider 
the sequence 
$$
d=(\frac{1}{2}, \frac{i}{2}, \frac{1+i}{2},0,1,i,0,1,i,0,1,i,\dots)\in\ell^\infty.  
$$
Then $d$ belongs to $\lm$ and the renormalized sum $s_X(d)$ vanishes.  
But there is no normal 
operator $N\in\mathcal N(X)$ whose diagonal relative to some orthonormal 
basis is $d$.  
\end{prop}

Before giving the proof, we recall that a doubly stochastic 
$n\times n$ matrix $A=(a_{ij})$ is said to be {\em orthostochastic} if 
there is a unitary $n\times n$ matrix $(u_{ij})$ such that 
$a_{ij}=|u_{ij}|^2$, $1\leq i,j\leq n$.  We will make use 
of the following known example of a doubly 
stochastic matrix that is not orthostochastic \cite{hornDiag}.  

\begin{lem}\label{cvLem1}
The $3\times 3$ matrix 
$$
(a_{ij})=
\begin{pmatrix}
1/2&0&1/2\\
0&1/2&1/2\\
1/2&1/2&0
\end{pmatrix}
$$
is not orthostochastic.  
\end{lem}

\begin{proof}  Indeed, 
if there were a unitary $3\times 3$ matrix $U=(u_{ij})$ such that 
$|u_{ij}|^2=a_{ij}$ for all $ij$, then $U$ must have the form 
$$
U=\frac{1}{\sqrt{2}}
\begin{pmatrix}
a&0&b\\
0&c&d\\
e&f&0
\end{pmatrix}
$$
with complex entries satisfying $|a|=|b|=|c|=|d|=|e|=|f|=1$.  
But the 
rows of such a matrix cannot be mutually orthogonal.
\end{proof}

\begin{proof}[Proof of Proposition \ref{cvProp1}]  Straightforward 
computations (that we omit) show that for the set $X=\{0,1,i\}$, 
the group $K_X$ and the obstruction group $\Gamma_X$ are 
given by 
$$
K_X=\mathbb Z+\mathbb Z\cdot i,\qquad \Gamma_X=\mathbb C/K_X\cong 
\mathbb R/\mathbb Z\oplus \mathbb R/\mathbb Z\cong \mathbb T^2.  
$$
Let $x$ be the sequence 
$$
x=(0,0,0,0,1,i,0,1,i,0,1,i,\dots).  
$$
Obviously $x_n\in X$ for every $n=1,2,\dots$, $d_n=x_n$ except 
for $n=1,2,3$,  and 
$$
\sum_{n=1}^\infty d_n-x_n=\frac{1}{2}+\frac{i}{2}+\frac{1+i}{2}=1+i\in 
\mathbb Z+\mathbb Z\cdot i=K_X.  
$$
Hence $d\in \lm$ and $s_X(d)=0$.  

Every operator $N\in\mathcal N(X)$ has the form 
$N=P+iQ$, where $P,Q$ are mutually orthogonal infinite rank 
projections such that $\mathbf 1-(P+Q)$ has infinite rank.  
Assuming that there is such an operator $N$ whose matrix 
relative to some orthonormal basis $e_1, e_2,\dots$ has diagonal 
$d=(d_1,d_2,\dots)$, we argue to a contradiction as follows.  
Let $p,q\in\ell^\infty$ be the real and imaginary parts of the 
sequence $d$
\begin{align*}
p &= (\frac{1}{2}, 0, \frac{1}{2},0,1,0,0,1,0,0,1,0,\dots)  \\
q & = (0, \frac{1}{2}, \frac{1}{2},0,0,1,0,0,1,0,0,1,\dots).  
\end{align*}
Since $d$ is the diagonal of $P+iQ$, 
one may equate real and imaginary parts to obtain 
$p_n=\langle Pe_n,e_n\rangle$ and 
$q_n=\langle Qe_n,e_n\rangle$, $n=1,2,\dots$.  

Now for $n\geq 4$, both $p_n$ and $q_n$ are $\{0,1\}$-valued.  Since 
$P$ and $Q$ are projections, it follows that 
$$
Pe_n=p_ne_n,\quad Qe_n=q_ne_n,\qquad n\geq 4,  
$$
and in particular, both $P$ and $Q$ leave the closed 
linear span $[e_4, e_5, e_6,\dots]$ invariant.  Hence 
they leave its orthocomplement $[e_1, e_2,e_3]$ invariant as well.  Let 
$P_0$, $Q_0$ be the restrictions of $P$, $Q$, respectively,  to 
$H_0=[e_1, e_2,e_3]$.  $P_0$ and $Q_0$ are mutually orthogonal
projections, 
and the diagonals of their matrices relative to the orthonormal 
basis $e_1, e_2, e_3$ are respectively 
$$
(\frac{1}{2}, 0, \frac{1}{2}), \quad{\rm and\ } (0,\frac{1}{2},\frac{1}{2}). 
$$
Each of these two diagonals has sum $1$, hence $P_0$ and 
$Q_0$ are one-dimensional.  Moreover,  $R_0=\mathbf 1_{H_0}-(P_0+Q_0)$ 
is a one-dimensional projection in $\mathcal B(H_0)$ 
whose diagonal relative to the basis 
$e_1, e_2, e_3$ is 
$$
(\frac{1}{2}, \frac{1}{2},0).  
$$
Hence the $3\times 3$ matrix whose rows are the diagonals of the 
three projections $P_0,Q_0,R_0$ takes the form
$$
A=
\begin{pmatrix}
1/2&0&1/2\\
0&1/2&1/2\\
1/2&1/2&0
\end{pmatrix}
.
$$
If we now choose unit vectors $f_1$, $f_2$, $f_3$ so that 
$P_0=[f_1]$, $Q_0=[f_2]$ and $R_0=[f_3]$, we find that $f_1, f_2, f_3$ is a second 
orthonormal basis for $H_0$, and
$$
A=
\begin{pmatrix}
|\langle e_1,f_1\rangle|^2&|\langle e_2,f_1\rangle|^2&|\langle e_3,f_1\rangle|^2\\
|\langle e_1,f_2\rangle|^2&|\langle e_2,f_2\rangle|^2&|\langle e_3,f_2\rangle|^2\\
|\langle e_1,f_3\rangle|^2&|\langle e_2,f_3\rangle|^2&|\langle e_3,f_3\rangle|^2
\end{pmatrix}
=
(|u_{ij}|^2)
, 
$$
where $(u_{ij})$ is a unitary $3\times 3$ matrix.  
This contradicts Lemma \ref{cvLem1}.  
\end{proof}

Proposition \ref{cvProp1} shows that the necessary condition $s_X(d)=0$ 
is not sufficient for a sequence 
$d\in \lm$ to be the diagonal of an operator in $\mathcal N(X)$ 
when $X$ contains more than two points.  Moreover, the precise  
nature of the remaining obstructions when $X$ consists of three 
non-colinear points remains mysterious.

\bibliographystyle{alpha}

\newcommand{\noopsort}[1]{} \newcommand{\printfirst}[2]{#1}
  \newcommand{\singleletter}[1]{#1} \newcommand{\switchargs}[2]{#2#1}


\begin{thebibliography}{Kad02b}

\bibitem[AK06]{arvKad1}
W.~Arveson and R.~Kadison.
\newblock Diagonals of self-adjoint operators.
\newblock In D.~R.~Larson D.~Han, P. E. T.~Jorgensen, editor, {\em Operator
  theory, operator algebras and applications}, Contemp Math. Amer. Math. Soc.,
  2006.
\newblock ar{X}iv:math.OA/0508482 v2.

\bibitem[Ati82]{atiyHorn}
M.~Atiyah.
\newblock Convexity and commuting {H}amiltonians.
\newblock {\em Bull. L.M.S.}, 14(1):1--15, 1982.

\bibitem[Eff89]{EdIntell}
E.~G. Effros.
\newblock Why the circle is connected.
\newblock {\em Math. Intelligenceer}, 11(1), 1989.

\bibitem[GS82]{GSI}
V.~Guillemin and S.~Sternberg.
\newblock Convexity properties of the moment mapping.
\newblock {\em Invent. Math.}, 67(3):491--513, 1982.

\bibitem[GS84]{GSII}
V.~Guillemin and S.~Sternberg.
\newblock Convexity properties of the moment mapping {II}.
\newblock {\em Invent. Math.}, 77(3):533--546, 1984.

\bibitem[Hor54]{hornDiag}
A.~Horn.
\newblock Doubly stochastic matrices and the diagonal of a rotation matrix.
\newblock {\em Amer. J. Math.}, 76(3):620--630, July 1954.

\bibitem[Kad02a]{kadPnasI}
R.~Kadison.
\newblock The {P}ythagorean theorem {I}: the finite case.
\newblock {\em Proc. N.A.S. (USA)}, 99(7):4178--4184, 2002.

\bibitem[Kad02b]{kadPnasII}
R.~Kadison.
\newblock The {P}ythagorean theorem {II}: the infinite discrete case.
\newblock {\em Proc. N.A.S. (USA)}, 99(8):5217--5222, 2002.

\bibitem[Kos73]{KostantSHthm}
B.~Kostant.
\newblock On convexity, the {W}eyl group and the {I}wasawa decomposition.
\newblock {\em Ann. Scient. {\'E}c. Norm. Sup.}, 6:413--455, 1973.

\bibitem[Neu99]{aNeumann1}
A.~Neumann.
\newblock An infinite-dimensional generalization of the {S}chur-{H}orn
  convexity theorem.
\newblock {\em Jour. Funct. Anal.}, 161(2):418--451, 1999.

\bibitem[Sch23]{schurDiag}
I.~Schur.
\newblock {\"U}ber eine klasse von mittlebildungen mit anwendungen auf der
  determinantentheorie.
\newblock {\em Sitzungsber. Berliner Mat. Ges.}, 22:9--29, 1923.

\end{thebibliography}

\end{document}